\documentclass[11pt,draft]{article}
\usepackage{amsmath}
\usepackage{amssymb}
\usepackage{latexsym}
\usepackage{theorem}
\usepackage{mathrsfs}
\usepackage[all]{xy}

\newcommand{\smashprod}{\wedge}

\newcommand{\zz}{{ \mathbb{Z} }}

\newcommand{\qq}{{ \mathbb{Q} }}

\newcommand{\dscr}{{ \mathscr{D} }}
\newcommand{\cscr}{{ \mathscr{C} }}
\newcommand{\ccal}{{ \mathcal{C} }}
\newcommand{\dcal}{{ \mathcal{D} }}
\newcommand{\ecal}{{ \mathcal{E} }}
\newcommand{\mcal}{{ \mathcal{M} }}
\newcommand{\gcal}{{ \mathcal{G} }}

\newcommand{\sscr}{{ \mathscr{S} }}

\newcommand{\tscr}{{ \mathscr{T} }}

\newcommand{\fscr}{{ \mathscr{F} }}

\newcommand{\fcal}{{ \mathcal{F} }}

\newcommand{\co}{{ : }}
\newcommand{\GIS}{{ G \mathscr{I} \mathscr{S} }}

\newcommand{\fibrep}{{ \widehat{f} }}
\newcommand{\cofrep}{{ \widehat{c} }}

\DeclareMathOperator{\id}{Id}
\DeclareMathOperator{\h}{H}
\DeclareMathOperator{\ho}{Ho}

\DeclareMathOperator{\leftmod}{--mod}

\DeclareMathOperator{\homSSS}{Hom}
\renewcommand{\hom}{\homSSS}
\DeclareMathOperator{\LimSSS}{lim}
\renewcommand{\lim}{\LimSSS}

\DeclareMathOperator{\SSET}{{sSet }}

\newtheorem{thm}{Theorem}[section]
\newtheorem{prop}[thm]{Proposition}
\newtheorem{cor}[thm]{Corollary}
\newtheorem{lem}[thm]{Lemma}

\newtheorem{defn}[thm]{Definition}

\newenvironment{pf}[1][Proof]{\vskip -0.25cm \textbf{#1} }
{\hfill \rule{0.5em}{0.5em}}

\begin{document}

\title{Splitting Monoidal Stable Model Categories}

\author{D. Barnes}

\maketitle

\begin{abstract}
\noindent If $\ccal$ is a stable model category
with a monoidal product then 
the set of homotopy classes of self-maps of the unit 
forms a commutative ring, $[S, S]^{\ccal}$. An idempotent $e$ of this
ring will split the homotopy category:
$[X, Y]^{\ccal} \cong e[X, Y]^{\ccal} \oplus (1-e)[X, Y]^{\ccal}$.
We prove that provided the localised model structures exist, this
splitting of the homotopy category comes from 
a splitting of the model category, that is, 
$\ccal$ is Quillen equivalent to $L_{e S} \ccal \times L_{(1-e)S} \ccal$
and $[X,Y]^{L_{eS} \ccal} \cong e [X,Y]^\ccal$.
This Quillen equivalence is strong monoidal and 
is symmetric when the monoidal product of $\ccal$ is. 
\end{abstract}

\section{Introduction}
Let $R$ be a commutative ring with an idempotent $e$, so $e \cdot e = e$, 
then there is an equivalence of categories 
$R \leftmod \overrightarrow{\longleftarrow} eR \leftmod \times (1-e) R \leftmod$
and for any $R$-module $M$ a natural isomorphism 
$M \cong eM \oplus (1-e)M$. 
This result can be useful since 
in general it is easier to study the categories 
$eR \leftmod$ and $(1-e) R \leftmod$ separately. 
We want to find some generalisation of this result to 
model categories. Our initial example is an additive and monoidal category, 
so we look for a class of monoidal model categories whose 
homotopy category is additive. The collection of 
monoidal stable model categories is such a class. 

A pointed model category $\ccal$ comes with a natural 
adjunction $(\Sigma, \Omega)$ on $\ho \ccal$. When this adjunction is an equivalence
we say that $\ccal$ is stable.
The homotopy category of a stable model category is naturally a 
triangulated category (hence additive), see \cite[Chapter 7]{hov99}.
We are interested in monoidal stable model categories:
those stable model categories which are also
monoidal model categories (\cite[Section 6.6]{hov99}). 
Thus $\ccal$ has a closed monoidal product $(\smashprod, \hom)$
with unit $S$ which is compatible with the model structure in the 
sense that the pushout product axiom holds. 
We write $[X,Y]^\ccal$ for the set of maps in the homotopy category of $\ccal$,
this is a group since $X$ is equivalent to $\Omega^2 \Sigma^2 X$. 
It is then an simple task to prove that 
$[S,S]^\ccal$ is a commutative ring (Lemma \ref{lem:unitring}). 

For any $X, Y \in \ccal$, 
$[X,Y]^\ccal$ is a $[S,S]^\ccal$-module via the smash product. 
Hence, for any idempotent $e \in [S,S]^\ccal$, 
we have an isomorphism which is natural in $X$ and $Y$:
$[X, Y]^{\ccal} \cong e[X, Y]^{\ccal} \oplus (1-e)[X, Y]^{\ccal}$.
Define $e \ho \ccal$ to be that category with the same class of objects
as $\ho \ccal$ and with morphisms given by $e[X, Y]^{\ccal}$.
Then, as with the case of $R$-modules above, we have an equivalence of categories
$\ho \ccal \overrightarrow{\longleftarrow} e \ho \ccal \times (1-e) \ho \ccal$.

We want to understand this splitting in terms of 
the model category $\ccal$. 
We assume that for any cofibrant object $E \in \ccal$
there is a new model structure on the category $\ccal$, written
$L_E \ccal$, with the same cofibrations as $\ccal$
and weak equivalences those maps $f$ such that
$\id_E \smashprod f$ is a weak equivalence of $\ccal$. 
The model structure $L_E \ccal$ is called the Bousfield localisation of
$\ccal$ at $E$ and there is a left Quillen functor
$\id : \ccal \to L_E \ccal$. 

For $e$ an idempotent of $[S,S]^\ccal$,
we are interested in localising at the objects 
$eS$ and $(1-e)S$. These are constructed in terms
of homotopy colimits and $S$ is weakly equivalent to
$eS \coprod (1-e)S$. 
Our main result, Theorem \ref{thm:generalsplitting}, is that the adjunction
$$\Delta :  \ccal \overrightarrow{\longleftarrow} L_{eS} \ccal \times L_{(1-e)S} \ccal : \prod$$
is a Quillen equivalence.
Furthermore $[X,Y]^{L_{eS} \ccal} \cong e [X,Y]^\ccal$,
so that this Quillen equivalence induces the splitting of 
$\ho \ccal$. 

Note that there is a non-trivial idempotent $e \in [S,S]^\ccal$
if and only if there is a non-trivial splitting of the
homotopy category. The splitting theorem proves that if there
is such an idempotent, then there is a splitting of model 
categories. Corollary \ref{cor:possiblesplittings}
demonstrates that if one has a splitting at the model category level
(into $L_E \ccal$ and $L_F \ccal$)
then the idempotent this defines ($e$) returns the splitting
at the model category level: $L_{eS} \ccal = L_E \ccal$
and $L_{(1-e)S} \ccal = L_F \ccal$. Hence, the notions:
a splitting of $[S,S]^\ccal$, a splitting of 
$\ho \ccal$ and a splitting of the model category 
$\ccal$, are all equivalent.

Our motivation for this splitting result came from 
studying rational equivariant spectra for compact Lie groups $G$. 
The ring of self-maps of the unit in the homotopy 
category of rational $G$-spectra,
$[S,S]^G_\qq$, is naturally isomorphic to the
rational Burnside ring. We have a good understanding of
idempotents in this ring via tom-Dieck's isomorphism, 
see Lemma \ref{lem:tdisom}. If a non-trivial idempotent
exists, then we can use it to split the category and obtain two pieces
which are possibly easier to study. 
We construct a model category of rational equivariant spectra
in Section \ref{sec:rateqspec}, we then give two 
examples of this splitting result
taken from \cite{barnes}.
Corollary \ref{cor:finG} considers the case of a finite group
and at the homotopy level recovers the splitting result of 
\cite[Appendix A]{gremay95}.
The second example is Lemma \ref{lem:Geidemfamily}
and in the case of $O(2)$ the idempotent constructed
is non-trivial and gives the homotopy 
level splitting of \cite{gre98a}.

Since we are working in a monoidal context and the splitting result is 
a strong monoidal adjunction, we can give 
two further examples: the case of modules
over a ring spectrum (Proposition \ref{prop:splitrmod})
and $R$-$R$-bimodules for a ring spectrum $R$
(Proposition \ref{prop:splitbimod}). After these examples we return to
our motivating case of rational $G$-spectra and give a model structure
for rational $G$-spectra in terms of modules over a commutative ring spectrum.

We also feel that we should mention \cite{ss03stabmodcat}. In this paper
the authors assume that one has a stable model category with a set of compact 
generators and conclude that such a category is Quillen equivalent to 
the category of right modules over a ring spectrum with many objects 
(that is, right modules over a category enriched over symmetric spectra). 
Consider a symmetric monoidal category $\ccal$ with a set of compact generators $\gcal$
such that there is an idempotent $e \in [S,S]^\ccal$, we can relate our 
splitting result to the work of the above-mentioned paper as follows. 
We have two new sets of compact objects 
$e \gcal = \{e G | G \in \gcal  \}$ and $(1-e) \gcal$,
their union is a set of generators for $\ccal$. 
We can construct a ring spectrum with many objects from 
$e \gcal$, call this $\ecal(e \gcal)$. 
The homotopy category of right modules over $\ecal(e \gcal)$
is equivalent to $e \ho \ccal$ and similarly the homotopy
category of right modules over $\ecal((1-e) \gcal)$
is equivalent to $(1-e) \ho \ccal$. All of our examples 
(see Sections \ref{sec:rateqspec} - \ref{sec:modbimod})
have a set of compact generators. 

\paragraph*{Acknowledgments} 

This work is a development of material from my PhD thesis, 
supervised by John Greenlees. I would to thank him for all the 
help and advice he has given me. 
Stating the splitting result in terms of stable model categories
was a suggestion made by both Neil Strickland and Stefan Schwede
and the bimodule example was the idea of Stefan Schwede, they both deserve
a great deal of gratitude for carefully reading earlier versions of this work.

\section{Stable Model Categories}

We introduce the notion of a stable model category, 
prove that if $\ccal$ is a monoidal stable model category
then $[S,S]^\ccal$ is a commutative ring and prove some basic results
about idempotents of $[S,S]^\ccal$. 

A pointed model category $\ccal$ comes with a natural action of 
$\ho \SSET_*$ (the homotopy category of pointed simplicial sets) 
on $\ho \ccal$, see \cite[Chapter 6]{hov99}
or \cite[Section I.2]{quil67}. In particular
for $X \in \ccal$ we have $\Sigma X := S^1 \smashprod^L X$ and
$\Omega X := R \hom_* (S^1, X)$, these define the 
suspension and loop adjunction $(\Sigma, \Omega)$
on $\ho \ccal$.
When this adjunction is an equivalence
we say that $\ccal$ is \textbf{stable}, see \cite[Chapter 7]{hov99}.
Following that chapter we see that $\ho \ccal$ is a triangulated 
category in the classical sense (see \cite{del77}) and that 
cofibre and fibre sequences agree (up to signs). 

Let $\ccal$ be a monoidal stable model category,
we let $\cofrep$ and $\fibrep$ denote cofibrant and fibrant replacement in $\ccal$.
For any collection of objects $\{Y_i \}_{i \in I}$ in $\ccal$, 
there is a natural map $\coprod Y_i \to \prod Y_i$. 
In a triangulated category finite coproducts and
finite products coincide, thus when 
$I$ is a finite set we have a weak equivalence
$\coprod_{i \in I} \cofrep Y_i \to \prod_{i \in I} \fibrep Y_i$.

\begin{lem}\label{lem:unitring}
The set $[S,S]^\ccal$ is a commutative ring. 
\end{lem}
\begin{pf}
The homotopy category of a stable model category is
additive \cite[Lemma 7.1.2]{hov99}. Thus 
$[S,S]^\ccal$ is an abelian group and
this addition is compatible with composition of maps
$\circ : [S,S]^\ccal \otimes_\zz [S,S]^\ccal \to [S,S]^\ccal$. 
There is also a smash product operation
$\smashprod : [S,S]^\ccal \otimes_\zz [S,S]^\ccal \to [S,S]^\ccal$.
The operations $\circ$ and $\smashprod$ satisfy the following interchange law.
Let $a$, $b$, $c$ and $d$ be elements of $[S,S]^\ccal$, then
$(a \circ b) \smashprod (c \circ d) =
(a \smashprod c) \circ (b \smashprod d) $
as elements of $[S \smashprod S, S\smashprod S]^\ccal$
and the unit of each operation is the identity map of $S$. 
Hence, by the well-known argument below, 
the two operations $\circ$ and $\smashprod$
are equal and commutative. So composition defines a 
commutative ring structure on the group
$[S,S]^\ccal$.

Consider any set $A$, with two binary operations
$\smashprod$, $\circ$ which satisfies the 
above interchange law. Assume there is an element 
$e \in A$ which acts as a both a left and right identity
for $\smashprod$ and $\circ$. 
Then 
$a \smashprod d = (a \circ e) \smashprod (e \circ d)
= a \circ d$
and $a \smashprod d = 
(e \circ a) \smashprod (d \circ e) =
d \circ a$. Hence the two operations are equal
and are commutative. 
\end{pf}

Note that the above does not assume that 
$\smashprod$ is a {\em symmetric} monoidal product. 
Consider a map in the homotopy category, $a \in [S,S]^\ccal$.
This can be represented by $a' \co \cofrep \fibrep S \to \cofrep \fibrep S$.
We can consider the homotopy colimit of the diagram
$\cofrep \fibrep S \overset{a'}{\to} 
\cofrep \fibrep S \overset{a'}{\to} \cofrep \fibrep S 
\overset{a'}{\to} \dots $
which we denote by $a S$. 
A different
choice of representative will give a weakly equivalent
homotopy colimit, so we must use a little care when writing
$a S$. The construction of the homotopy colimit $a S$
comes with a map $\cofrep S \to \cofrep \fibrep S \to a' \cofrep \fibrep S$.
For any $X \in \ccal$, 
we have the map $a' \smashprod \id_X \co \cofrep \fibrep S \smashprod X \to 
\cofrep \fibrep S \smashprod X$. We can then construct homotopy colimits as above
to create the object $aX$.
We use \cite[Proposition 7.3.2]{hov99}, to obtain an 
exact sequence:
$$0 \to \lim^1 [X, Y]^\ccal \to [a X,Y]^\ccal \to \lim [X,Y]^\ccal \to 0.$$
We are interested in $e S$ for $e$ an idempotent of $[S,S]^\ccal$. 
In such a case, the $\lim^1$-term is zero as 
the tower created by an idempotent satisfies the 
Mittag-Leffler condition (\cite[Definition 3.5.6]{weib}). Hence the above 
exact sequence reduces to an isomorphism 
$[e X,Y]^\ccal \to \lim [X,Y]^\ccal = e [X,Y]^\ccal$.

If $e$ is an idempotent so is $(\id_S - e)$, which we now write as 
$(1-e)$. Furthermore we have a canonical natural isomorphism 
$[X,Y]^\ccal \cong e[X,Y]^\ccal \oplus (1-e)[X,Y]^\ccal$ for any $X$ and $Y$.
Thus, there is a natural isomorphism in the homotopy category 
$X \to e X \prod (1-e) X$. 
We can write $\ho \ccal$ as the product
category $e \ho \ccal \times (1-e) \ho \ccal$, where 
$e \ho \ccal$ has the same objects as $\ho \ccal$
and $e \ho \ccal(X,Y) :=e[X,Y]^\ccal$.
We wish to pull this splitting back to the level of model
categories. 

\begin{lem}\label{lem:splitobjects}
For any object $X$ in $\ccal$
there is a natural weak equivalence
$\cofrep X \smashprod \cofrep S \to \fibrep e X \prod \fibrep (1-e) X$.
\end{lem}
\begin{pf}
We start with the maps
$\cofrep X \smashprod \cofrep S \to e X$
and $\cofrep X \smashprod \cofrep S \to (1-e) X$.
By taking fibrant replacements we obtain a map 
$\cofrep X \smashprod \cofrep S \to \fibrep e X \prod \fibrep (1-e) X$.
The following diagram commutes for any $Y \in \ccal$,
proving the result. 

$$\xymatrix{
[\fibrep eX \prod \fibrep (1-e) X, Y]^\ccal \ar[r] \ar[d]^\cong &
[X,Y]^\ccal \ar[dd]^\cong \\
[\fibrep eX \vee \fibrep (1-e) X, Y]^\ccal \ar[d]^\cong \\
[\fibrep eX , Y]^\ccal \oplus [\fibrep (1-e) X, Y]^\ccal \ar[r]^\cong &
e[X , Y]^\ccal \oplus (1-e) [X, Y]^\ccal }$$ 

\end{pf}

\section{Localisations}

We define the notion of a Bousfield localisation of a monoidal
model category and prove that when the localisation exists, 
the new model category shares many of the properties of the original
(left properness, the pushout product axiom and the monoid axiom).
We also consider Quillen pairs between localised categories.

Recall the following concepts of localisation.
\begin{defn}\label{def:genEstuff}
Let $E$ be a cofibrant object of the monoidal model category $\ccal$
and let $X$, $Y$ and $Z$ be objects of $\ccal$.
\begin{enumerate}
\item A map $f : X \to Y$ is an
$E$-\textbf{equivalence}\index{E-equivalence@$E$-equivalence}
if $\id_E \smashprod f : E \smashprod X \to E \smashprod Y$
is a weak equivalence.
\item $Z$ is $E$-\textbf{local}\index{E-local@$E$-local}
if $f^* : [Y,Z]^\ccal \to [X,Z]^\ccal$ is an isomorphism for all $E$-equivalences
$f : X \to Y$.
\item An $E$-\textbf{localisation}\index{E-localisation@$E$-localisation}
of $X$ is an $E$-equivalence $\lambda : X \to Y$ from $X$ to an $E$-local object $Y$.
\item $A$ is $E$-\textbf{acyclic}\index{E-acyclic@$E$-acyclic} 
if the map $* \to A$ is an $E$-equivalence.
\end{enumerate}
\end{defn}

The following is a standard result, see \cite[Theorems 3.2.13 and 3.2.14]{hir03}.
\begin{lem}\label{lem:genEequivElocal}
An $E$-equivalence between $E$-local objects is a weak equivalence.
\end{lem}

Consider the category $\ccal$ 
with a new set of weak equivalences:
the $E$-equivalences, while leaving the cofibrations unchanged.
If this defines a model structure we call this the Bousfield
localisation of $\ccal$ at $E$ and write it as $L_E \ccal$. 
The identity functor gives a strong monoidal Quillen pair (see definition below)
$$\id_\ccal : \ccal \overrightarrow{\longleftarrow} L_E \ccal : \id_\ccal.$$
This follows since the cofibrations are unchanged and 
if $f : X \to Y$ is an acyclic cofibration of $\ccal$
then $f \smashprod \id_E$ is also an acyclic cofibration.
Hence $f$ is a cofibration and an $E$-equivalence. 
We will write $\fibrep_{E}$ for fibrant replacement
in $L_E \ccal$.

\begin{defn}
A Quillen pair $L : \ccal \overrightarrow{\longleftarrow} \dcal : R$
between monoidal model categories is said to be a \textbf{strong monoidal adjunction}
if there is a natural isomorphism $L (X \otimes Y) \to LX \otimes LY$
and an isomorphism $L S_\ccal \to S_\dcal$. We require that these isomorphisms
satisfy the associativity and unital coherence conditions of
\cite[Definition 4.1.2]{hov99}. A strong monoidal adjunction $(L,R)$
is a \textbf{strong monoidal Quillen pair} if it is a Quillen adjunction
and if whenever $\cofrep S_\ccal \to S_\ccal$ is a cofibrant replacement
of $S_\ccal$, then the induced map
$L \cofrep S_\ccal \to L S_\ccal$ is a weak equivalence. 
\end{defn}

From now on we assume that for any cofibrant $E$
the $E$-equivalences and cofibrations define a model 
structure on $\ccal$, the $E$-local model structure.
In general we won't have a good description of the fibrations of 
$L_E \ccal$, however we do have the following lemma. This result is similar in nature to 
\cite[Proposition 3.4.1]{hir03}. 

\begin{lem}
An $E$-fibrant object is fibrant in $\ccal$ and $E$-local. 
If $X$ is $E$-local and fibrant in $\ccal$, then 
$X \to *$ has the right lifting property with respect to the
class of $E$-acyclic cofibrations between cofibrant objects.
\end{lem}

Note that in many cases a stronger result holds:
an object is $E$-fibrant if and only if it is fibrant in $\ccal$ and $E$-local. 
For example, this stronger result holds for EKMM spectra localised at an object $E$
by the fact that
the domains of the generating $E$-acyclic cofibrations
are cofibrant. \vskip 0.5cm

\begin{pf}
Let $A \to B$ be an acyclic cofibration, then 
this is also an $E$-equivalence. So for an $E$-fibrant object 
$Z$, the canonical map $Z \to *$ will have the right lifting 
property with respect to $A \to B$.
Let $f \co A \to B$ be an $E$-equivalence. We must prove that
$f^* \co [B,Z]^{\ccal} \to [A, Z]^{\ccal}$ is an isomorphism. 
But since $Z$ is $E$-fibrant the Quillen pair between $\ccal$
and $L_E \ccal$ gives an isomorphism
$[B,Z]^{\ccal} \cong [B,Z]^{L_E \ccal}$. This is natural in 
the first variable and the first statement follows.

Let $i : A \to B$ be an $E$-acyclic cofibration between cofibrant objects
and let $f : A \to X$ be any map of $\ccal$. 
Since $X$ is $E$-local, $i$ induces an isomorphism
$i^* :[B,X]^\ccal \to [A,X]^\ccal$. Choose $g : B \to X$
such that $g \circ i$ and $f$ are homotopic. We now apply the homotopy
extension property (see \cite[Page 1.7]{quil67}), choose a path object $X'$ for $X$
with a map $h : A \to X'$ such that $p_0 \circ h = g \circ i$
and $p_1 \circ h = f$. We thus have the following diagram
$$\xymatrix{
A \ar[r]^h \ar[d]_i & X' \ar[d]^{p_0} \\
B \ar[r]^g & X 
}$$
where $i$ is a cofibration and $p_0$ is a fibration
and a weak equivalence in $\ccal$.
Thus we have a lifting $H : B \to X'$ and the map
$p_1 \circ H$ is the solution to our original lifting problem.  
\end{pf}

If the $E$-local model structure
exists, then every weak equivalence is an $E$-equivalence.
Take a weak equivalence $f$, factor this into 
$g \circ h$ with $h$ a cofibration and a weak equivalence
and $h$ an acyclic $E$-fibration. 
Then since smashing with $E$ is a left Quillen functor, 
$\id_E \smashprod h$ is an acyclic cofibration.
By definition, $\id_E \smashprod g$ is a weak equivalence,
hence so is $\id_E \smashprod f$.
We also note that if $F$ and $E$ are cofibrant objects of $\ccal$
then the model categories $L_{F \smashprod E} \ccal$ and $L_E L_F \ccal$
are equal (they have the same weak equivalences and cofibrations).

Now we prove a straightforward result about Quillen functors
between localised categories and then turn to proving that
$L_E \ccal$ inherits many of the properties of the original model structure
on $\ccal$.

\begin{thm}\label{thm:locfuncs}
Take a Quillen adjunction between monoidal model categories with 
a strong monoidal left adjoint
$F : \cscr \overrightarrow{\longleftarrow} \dscr : G.$
Let $E$ be cofibrant in $\cscr$ and assume that all model categories
mentioned below exist. Then $(F,G)$ passes to a Quillen pair
$F : L_E \cscr \overrightarrow{\longleftarrow} L_{FE} \dscr : G.$
Furthermore, if $(F,G)$ form a Quillen equivalence,
then they pass to a Quillen equivalence of the localised categories.
\end{thm}
\begin{pf}
Since the cofibrations in $L_E \cscr$ and $L_{FE} \dscr$ are unchanged
$F$ preserves cofibrations. Now take an acyclic cofibration in $\cscr$ of the form
$\id_E \smashprod f : E \smashprod X \to E \smashprod Y$, applying $F$ and
using the strong monoidal condition we have a weak equivalence in $\dscr$:
$\id_{FE} \smashprod Ff : FE \smashprod FX \to FE \smashprod FY$.
Hence $F$ takes $E$-acyclic cofibrations to $FE$-acyclic cofibrations and we have a Quillen pair.

To prove the second statement we
show that $F$ reflects $E$-equivalences
between cofibrant objects and that $F \cofrep GX \to X$ is an $E$-equivalence
for all $X$ fibrant in $L_{FE} \dscr$. These conditions are an equivalent
definition of Quillen equivalence by \cite[Corollary 1.3.16(b)]{hov99}.
The first condition follows since strong monoidality
allows us to identify $F(\id_E \smashprod f)$
and $\id_{FE} \smashprod Ff $ for a map $f$ in $\cscr$
and $F$ reflects weak equivalences between cofibrant objects.
The second condition is equally simple:
we know that an $E$-fibrant object is fibrant
and that cofibrant replacement is unaffected
by Bousfield localisation. Hence
$F \cofrep GX \to X$ is a weak equivalence
and thus an $E$-equivalence. 
\end{pf}

\begin{prop}\index{Left proper}
If $\ccal$ is left proper so is $L_E \ccal$.
\end{prop}

\begin{prop}\label{prop:pushaxiom}
If $\ccal$ is symmetric monoidal, then
for two cofibrations, $f : U \to V$ and $g : W \to X$, the induced map
$$f \square g : V \smashprod  W \bigvee_{U \smashprod W} U \smashprod X  \to V \smashprod X$$
is a cofibration which is an $E$-acyclic cofibration if either $f$ or $g$ is.
If $X$ is a cofibrant object then the map $\cofrep S \smashprod X \to X$
is a weak equivalence.
\end{prop}
\begin{pf}
Since the cofibrations are unchanged by localisation,
we only need to check that the above map is an
$E$-equivalence when one of $f$ or $g$ is.
Assume that $f$ is an $E$-equivalence, then 
the map $\id_E \smashprod f : E \smashprod U
\to E \smashprod V$ is a weak equivalence
and a cofibration. Thus, since $E \smashprod (-)$ commutes with
pushouts the map
$$ E \smashprod (V \smashprod  W \bigvee_{U \smashprod W} U \smashprod X)
  \to E \smashprod (V \smashprod X) $$
is also a weak equivalence and a cofibration.
By symmetry, this also deals with the case when $g$ is an $E$-equivalence.
The unit condition is unaffected by localisation, so it holds
in the $E$-local model structure. 
\end{pf}

Thus, when $\ccal$ is symmetric, $L_E \ccal$ is a monoidal model category.
Now we consider the monoid axiom.

\begin{prop}\label{prop:Emonoid}
If $\ccal$ is symmetric monoidal and satisfies the monoid axiom, then so does 
$L_E \ccal$. 
\end{prop}
\begin{pf}
Let $i : A \to X$ be an acyclic $E$-cofibration, then for any object $Y$,
the map $\id_E \smashprod i \smashprod \id_Y$ 
is a weak equivalence. 
Moreover, transfinite compositions of pushouts of such maps are
weak equivalences by the monoid axiom for $\ccal$. 
Thus transfinite compositions of pushouts of maps of the form 
$i \smashprod \id_Y$ are $E$-equivalences. 
\end{pf}

\section{The Splitting}

We are now ready to prove our main result, 
Theorem \ref{thm:generalsplitting}. We conclude this
section with a converse to this result. 

Recall the definition of the
product model category\index{Product model category}
from \cite[Example 1.1.6]{hov99}.
Given model categories $M_1$ and $M_2$ we can put a model category structure on
$M_1 \times M_2$. A map $(f_1,f_2)$ is a cofibration, weak equivalence or fibration
if and only if $f_1$ is so in $M_1$ and $f_2$ is so in $M_2$. Similarly
a finite product of model categories has a model structure where a map
is a cofibration, weak equivalence or fibration
if and only if each of its factors is so.
If $M_1$ and $M_2$ both satisfy any of the following:
left properness, right properness, the pushout product axiom, the monoid axiom or
cofibrant generation, then so does $M_1 \times M_2$.

\begin{prop}\label{prop:genadjunct}
If $E$ and $F$ are cofibrant objects of $\ccal$ 
then there is a strong monoidal Quillen adjunction
$$\Delta :  \ccal \overrightarrow{\longleftarrow} L_E \ccal \times L_F \ccal : \prod.$$
\end{prop}

Let $\ccal$ be a stable monoidal model category with 
an idempotent $e \in [S,S]^\ccal$. Then we have 
a Quillen pair 
$$\Delta :  \ccal \overrightarrow{\longleftarrow} L_{eS} \ccal \times L_{(1-e)S} \ccal : \prod$$
and an equivalence of homotopy categories
$$\Delta :  \ho \ccal \overrightarrow{\longleftarrow} e \ho \ccal \times (1-e) \ho \ccal : \prod.$$
We now wish to prove that the Quillen pair induces this 
equivalence of homotopy categories.

\begin{lem}
Take an idempotent $e \in [S,S]^\ccal$, any pair of
objects $X$, $Y$ and an $eS$-local object $Z$.
Then there are natural
isomorphisms 
$$[X, Y]^{L_eS \ccal} \longrightarrow [X,\fibrep_{eS}Y]^\ccal, \quad 
[X,Z]^{\ccal} \longrightarrow e[X,Z]^{\ccal}.$$
\end{lem}
\begin{pf}
The first comes from the Quillen adjunction between 
$\ccal$ and $L_{eS} \ccal$. 
For the second we use the fact that the map $\cofrep X \to e X$
is an $eS$-equivalence to obtain isomorphisms 
$[X,Z]^{\ccal} \leftarrow [eX,Z]^{\ccal} \to e[X,Z]^{\ccal}$. 
\end{pf}

\begin{lem}\label{lem:cofibreseqce}
Let $e$ be an idempotent of $[S,S]^\ccal$.
Then the map $e \co eS \to eS$ is an isomorphism
in $\ho \ccal$. Hence $(1-e) \co eS \to eS$
is equal to the zero map in $\ho \ccal$ and so 
for any $X$ and $Y$ in $\ccal$, $(1-e)[X,eY]^\ccal=0$.
\end{lem}
\begin{pf}
Consider the map $e^* \co [eS,X]^\ccal \to [eS,X]^\ccal$, this is naturally 
isomorphic to $e^* \co e[S,X]^\ccal \to e[S,X]^\ccal$, which is an isomorphism.   
The second part follows since $(1-e) \circ e \in [S,S]^\ccal$
is equal to zero. 
\end{pf}

\begin{thm}\label{thm:generalsplitting}
Let $\ccal$ be a stable monoidal model category with 
an idempotent $e \in [S,S]^\ccal$. 
Assume that the model categories 
$L_{eS} \ccal$ and $L_{(1-e)S} \ccal$ exist, 
then the strong monoidal Quillen pair below  is a Quillen equivalence. 
$$\Delta :  \ccal \overrightarrow{\longleftarrow} 
L_{eS} \ccal \times L_{(1-e)S} \ccal : \prod$$
\end{thm}
\begin{pf}
The right adjoint detects all weak equivalences:
take $f \co A \to B$ in $L_{eS} \ccal$ and $g \co C \to D$ in $L_{(1-e)S} \ccal$. If 
$(f,g) \co A \prod C \to B \prod D$ is a weak equivalence
then $f$ and $g$ are weak equivalences since they are retracts of $(f,g)$. 
Hence $f$ is an $eS$-equivalence 
and $g$ is a $(1-e)S$-equivalence. 

Let $X$ be a cofibrant object of $\ccal$, 
we then have an $eS$-acyclic cofibration $X \to \fibrep_{eS} X$
and an $(1-e)S$-acyclic cofibration $X \to \fibrep_{(1-e)S} X$. We must prove that
$X \to \fibrep_{eS} X \prod \fibrep_{(1-e)S} X$ is a weak equivalence.
For any $A \in \ccal$ we have the following commutative diagram:
$$\xymatrix{
e[A,X]^\ccal \oplus (1-e)[A,X]^\ccal \ar[r] &
e[A,\fibrep_{eS} X]^\ccal \oplus (1-e)[A,\fibrep_{(1-e)S} X]^\ccal \\
[A, X]^\ccal \ar[r] \ar[u]_\cong &
[A, \fibrep_{eS} X \prod \fibrep_{(1-e)S} X]^\ccal \ar[u]_\cong
}$$
So we have reduced the problem to proving that 
$e[A,X]^\ccal \to e[A,\fibrep_{eS} X]^\ccal$ is an isomorphism. 
This follows from the commutative diagram below and Lemma \ref{lem:cofibreseqce}, 
which tells us that the terms $ e[A,(1-e)X]^\ccal$ and
 $e[A,(1-e)\fibrep_{eS}X]^\ccal$ are zero. 
$$\xymatrix@C+0.3cm{
e[A,X]^\ccal \ar[r] \ar[d]^\cong & 
e[A,\fibrep_{eS} X]^\ccal \ar[d]^\cong \\
e[A,eX]^\ccal \oplus e[A,(1-e)X]^\ccal \ar[r]^(0.45){\cong}  & 
e[A, e\fibrep_{eS} X]^\ccal \oplus e[A,(1-e)\fibrep_{eS}X]^\ccal
}$$ 
\end{pf}

A \textbf{finite orthogonal decomposition} of $\id_S$
is a collection of idempotents
$e_1, \dots, e_n$ which sum to the identity in $[S,S]^\ccal$ such that 
$e_i \circ e_j =0$ for $i \neq j$.
This result extends to give a strong monoidal
Quillen equivalence between $\ccal$ and $\prod_{i=1}^n L_{e_i S} \ccal$
whenever $e_1, \dots, e_n$ is a finite orthogonal decomposition of $\id_S$.

\begin{cor}\label{cor:possiblesplittings}
Consider a monoidal model category $\ccal$ which splits as a product 
$L_E \ccal \times L_F \ccal$, for cofibrant objects  $E$ and $F$.
Then there are orthogonal idempotents
$e_E$ and $e_F$ in $[S,S]^\ccal$ such that $e_E + e_F = \id_S$, 
$L_{e_E S} \ccal = L_E \ccal$ and $L_{e_F S} \ccal = L_F \ccal$.
\end{cor}
\begin{pf}
Using the isomorphism $[S,S]^{L_E \ccal} \oplus [S,S]^{L_F \ccal} \to [S,S]^\ccal$
define $e_E$ as the image of $\id_S \oplus \ 0 \in [S,S]^{L_E \ccal} \oplus [S,S]^{L_F \ccal}$
in $[S,S]^\ccal$. Similarly define $e_F$ as the image of $0 \oplus \id_S$.
Thus we have idempotents $e_E$ and $e_F$ in $[S,S]^\ccal$ such that 
$e_E + e_F = \id_S$ and $e_E \circ e_F = 0$. 
By construction, $e_E [X,Y]^\ccal \cong [X,Y]^{L_E \ccal}$ and by our work above
$e_E [X,Y]^\ccal \cong [X,Y]^{L_{e_E S} \ccal}$. 
From this it follows that the $e_E S$-equivalences are the 
$E$-equivalences and $L_{eS} \ccal =L_E \ccal$. 
\end{pf}

\section{Rational Equivariant Spectra}\label{sec:rateqspec}
Our motivating example for the splitting result is the category
of rational $G$-equivariant EKMM $S$-modules for a compact Lie group $G$.
Our first task is to define this category, for this we will need a 
rational sphere spectrum.
We work with $G \mcal$, the category of 
$G$-equivariant EKMM $S$-modules from  \cite{mm02}. 
One could work with $G$-equivariant orthogonal spectra
and perform analogous constructions there and obtain
equivalent results for that category. 
In particular the two categories of equivariant
spectra we have mentioned are monoidally Quillen equivalent.  

We will construct $\qq$ as a group and translate this into spectra.
Take a free resolution of $\qq$ as an abelian group, $0 \to R
\overset{f}{\to} F \to \qq \to 0$, where $F= \oplus_{q \in \qq} \zz$. 
Since a free abelian group is
a direct sum of copies of $\zz$ we can rewrite this short exact
sequence as $0 \to \bigoplus_i \zz \overset{f}{\to} \bigoplus_j
\zz \to \qq \to 0$. Since $\qq$ is flat, the sequence $0 \to
\bigoplus_i M \overset{f \otimes \id}{\to} \bigoplus_j M \to \qq \otimes M \to
0$ is exact for any abelian group $M$.
Hence for each subgroup $H$ of $G$, we have an injective  map
(which we also denote as $f$)
$\bigoplus_i A(H) \overset{f \otimes \id}{\to} \bigoplus_j A(H)$
and $\bigoplus_j A(H) / \bigoplus_i A(H) \cong A(H) \otimes \qq$.
For $H$, a subgroup of $G$,
$$[\bigvee_i S, \bigvee_j S]^H
\cong \hom_{A(H)} \Big( \bigoplus_{i} A(H),\bigoplus_{j} A(H) \Big).$$
Thus we can choose
$g : \bigvee_i \cofrep S \to \bigvee_j \cofrep S$,
a representative for the homotopy class corresponding to $f$.

Let $I$ be the unit interval with basepoint $0$, 
there is a cofibration of spaces $S^0 \to I$
which sends the non-basepoint point of $S^0$ to $1 \in I$. 
If $X$ is a cofibrant $G$-spectrum then 
$X \cong X \smashprod S^0 \to X \smashprod I$ 
is a cofibration since $G$-spectra are enriched over spaces
(see \cite[Chapter III, Definition 1.14]{mm02} and 
\cite[Lemma 4.2.2]{hov99}).
For a map $f : X \to Y$, the cofibre
of $f$, $C_f$, is the pushout of the diagram 
$X \smashprod I \leftarrow X \overset{f}{\to} Y$. 
If $X$ is cofibrant then the map 
$Y \to C_f$ is a cofibration, hence if $X$ and $Y$
are cofibrant, so is $Cf$. 

\begin{defn}\label{def:ratsphere}
For the map $g$ as constructed above, the cofibre of $g$ is the
\textbf{rational sphere spectrum}\index{Rational
sphere spectrum}\index{S rational@${S^0_\mcal{\qq}}$}
and we have a cofibre sequence
$$\bigvee_i \cofrep S \overset{g}{\longrightarrow}
\bigvee_j \cofrep S \longrightarrow S^0_\mcal {\qq}. $$
\end{defn}
A different choice of representative for the homotopy class $[g]$
will induce a weak equivalence between the cofibres,
and hence (up to weak equivalence) $S^0_\mcal{\qq} $ is independent of
this choice of representative.
Note that there is an inclusion $\alpha : \cofrep S \to \bigvee_j \cofrep S$
which sends $\cofrep S$ to the term of $\bigvee_j \cofrep S$
corresponding to $1 \in \qq$. 

\begin{prop}\label{prop:rathomgps}
Let $X$ be a $G$-spectrum, then
for any subgroup $H$ of $G$ the map 
$(\id_X \smashprod \alpha)_* : \pi_*^H(X) \to \pi_*^H(X \smashprod S^0_\mcal{\qq} )$
induces an isomorphism 
$\pi_*^H(X) \otimes \qq \to \pi_*^H(X \smashprod S^0_\mcal{\qq} )$.
\end{prop}
\begin{pf}
Using the cofibre sequence which defines
$S^0{\qq}$ we have the following collection of isomorphic
long exact sequences of homotopy groups
$$\begin{array}{ccccccc}
\dots \longrightarrow &
\pi_n^H(X \smashprod \bigvee_i \cofrep S) & \overset{(\id \smashprod g)_*}{\longrightarrow} &
\pi_n^H(X \smashprod \bigvee_j \cofrep S) & \longrightarrow &
\pi_n^H(X \smashprod S^0_\mcal{\qq})    & \longrightarrow
\dots \\
\dots \longrightarrow &
\pi_n^H(\bigvee_i X) & \overset{(\id \smashprod g)_*}{\longrightarrow} &
\pi_n^H(\bigvee_j X) & \longrightarrow &
\pi_n^H(X \smashprod S^0_\mcal{\qq}) & \longrightarrow
\dots \\
\dots \longrightarrow &
\bigoplus_i \pi_n^H( X) & \overset{g \otimes \id}{\longrightarrow} &
\bigoplus_j \pi_n^H( X) & \longrightarrow &
\pi_n^H(X \smashprod S^0_\mcal{\qq}) & \longrightarrow
\dots \\
\dots \longrightarrow &
\bigoplus_i \zz \bigotimes \pi_n^H( X)   & \overset{g \otimes \id}{\longrightarrow} &
\bigoplus_j \zz \bigotimes \pi_n^H( X)   & \longrightarrow &
\pi_n^H(X \smashprod S^0_\mcal{\qq}) & \longrightarrow
\dots
\end{array} $$
Since the map $g \otimes \id : (\bigoplus_i \zz) \otimes \pi_n^H( X) \to
(\bigoplus_j \zz) \otimes \pi_n^H( X)$ is injective for all $n$,
this long exact sequence splits into short exact sequences
and the result follows. 
\end{pf}

There are many other methods for constructing a rational sphere spectrum,
these will all be weakly equivalent to $S^0_\mcal \qq$ as we prove below.
One obvious alternative is to construct a homotopy colimit of the diagram
$\cofrep S \overset{2}{\to} \cofrep S 
\overset{3}{\to} \cofrep  S \overset{4}{\to} \dots$, 
call this object $R_\qq$. It follows that the map
$\pi_*^H(\cofrep X) \to \pi_*^H(R_\qq \smashprod \cofrep X)$
induced by $\cofrep S \to R_\qq$ gives an isomorphism 
$\pi_*^H(\cofrep X) \otimes \qq \to \pi_*^H(R_\qq \smashprod \cofrep X)$.
We prove in Lemma \ref{lem:univpropsq} that if you have any 
rationalisation of the sphere -- a rational equivalence $f \co S \to X$
where $X$ is a spectrum with $\pi_*^H(X)$ rational
for all $n$ and $H$, then $S^0_\mcal \qq$
and $X$ are weakly equivalent. 

The result below is \cite[Chapter IV, Theorem 6.3]{mm02},
the proof of which is an adaptation
of the material in \cite[chapter VIII]{EKMM97}.
\begin{thm}\label{thm:GSlocal}
Let $E$ be a cofibrant spectrum or a cofibrant based $G$-space. Then $G \mcal $
has an
$E$-model structure\index{E-model structure@$E$-model structure}
whose weak equivalences are the $E$-equivalences and whose $E$-cofibrations are the
cofibrations of $G \mcal$.
The $E$-fibrant objects are precisely the $E$-local
objects and $E$-fibrant approximation
constructs a Bousfield localisation $f_X : X \to \fibrep_E X$ of $X$ at $E$.
The notation for $E$-model structure on the underlying category of $G \mcal$ is
$L_E G\mcal $ or $G \mcal_E$\index{G MM@$G \mcal_E$}.
\end{thm}

The categories $L_E G \mcal$ are cofibrantly generated model categories,
this is implied by the proof of
\cite[Chapter VIII, Theorem 1.1]{EKMM97}.
Let $c$ be a fixed infinite cardinal that
is at least the cardinality of
$E^*(S)$.
Then define $\tscr$, a test set for $E$-fibrations,
to consist of all inclusions of cell complexes
$X \to Y$ such that the cardinality of the set of cells
of $Y$ is less than or equal to $c$.
Hence the domains of these maps are $\kappa$-small
where $\kappa$ is the least cardinal greater than $c$.
Thus if we let $I$ be the set of generating cofibrations for $G \mcal$, then
we can take $I$ and $\tscr$ as sets of generating cofibrations and 
generating acyclic cofibrations
for $L_E G \mcal$.

\begin{lem}\label{lem:ratequivs}
For a map $g : X \to Y$ the following are equivalent:
\begin{enumerate}
\item $g : X \to Y $ is an $S^0_\mcal \qq$-equivalence.
\item $g_*^H : \pi_*(X^H) \otimes \qq \to
\pi_*(Y^H) \otimes \qq$ is an isomorphism for all $H$.
\item \vskip -0.0cm $g_*^H : \h_*(X^H; \qq)
\to \h_*(Y^H; \qq)$ is an isomorphism for all $H$.
\end{enumerate}
\end{lem}
\begin{pf}
We have shown in Proposition \ref{prop:rathomgps} that the first
two conditions are equivalent.
The last two statements are equivalent since the Hurewicz map
induces an isomorphism $\pi_*(A) \otimes \qq \to H_*(A;\qq)$
for any non-equivariant spectrum $A$.
\end{pf}

\begin{defn}
The model category of rational $G$-spectra is defined to be 
$L_{{S^0_\mcal \qq}} G \mcal$, which we write as $G \mcal_\qq$.
Since the $S^0_\mcal \qq$-equivalences
are precisely the rational homotopy isomorphisms, 
we call the $S^0_\mcal \qq$-equivalences
\textbf{rational equivalences} or \textbf{$\pi_*^\qq$-isomorphisms}.
The set of rational homotopy classes of maps from $X$ to $Y$ will be written $[X,Y]^G_\qq$
and we will write $\fibrep_\qq$ for fibrant replacement in the localised category.
\end{defn}
The lemma above proves that our model structure
is independent of our choice of rational
sphere spectrum.
We now prove that $G \mcal_\qq$ is a right proper model category,
for which we need the following.

\begin{lem}\label{lem:ratLES}
For any map $f : X \to Y$ of $G$-prespectra and any $H \subset G$, there are natural
long exact sequences
$$\xymatrix@C-0.42cm@R-0.6cm{
\dots \ar[r] &
\pi_q^H(Ff)     \otimes \qq \ar[r] &
\pi_q^H(X)      \otimes \qq \ar[r] &
\pi_q^H(Y)      \otimes \qq \ar[r] &
\pi_{q-1}^H(Ff) \otimes \qq \ar[r] &
\dots, \\
\dots \ar[r] &
\pi_q^H(X)      \otimes \qq \ar[r] &
\pi_q^H(Y)      \otimes \qq \ar[r] &
\pi_q^H(Cf)     \otimes \qq \ar[r] &
\pi_{q-1}^H(X)  \otimes \qq \ar[r] &
\dots }$$
and the natural map $\nu : Ff \to \Omega Cf$ is a $\pi_*$-isomorphism.
\end{lem}
\begin{pf}
By \cite[Chapter IV, Remark 2.8]{mm02}, we have long exact 
sequences as above, but without needing to tensor with $\qq$.
Since $\qq$ is flat, tensoring with it preserves exact sequences,
hence the result follows. 
\end{pf}

\begin{lem}\label{lem:rightproperrational}
The category $G \mcal_\qq$
is right proper.
\end{lem}
\begin{pf}
Following the proof of \cite[Lemma 9.10]{mmss01}
one shows that a stronger statement holds:
in a pullback diagram as below, if $\beta$
is a level wise fibration of $G$-spaces then $r$ is a $\pi_*^\qq$-isomorphism.
$$\xymatrix@!C{
W \ar[r]^{\delta} \ar[d]_{r} & X \ar[d]^{\sim_\qq} \\
Y  \ar[r]_\beta & Z
\ar@{}[ul]|\lrcorner|(0.52){\cdot \hskip 2.5pt } }$$
The only point of difference is that in the last step
of the proof one needs to use the long exact sequence
of \emph{rational} homotopy groups of a fibration. 
\end{pf}

Since is our localisation is of a particularly nice form, 
we are able to give the following interpretation 
of maps in $\ho G \mcal_\qq$. 

\begin{thm}\label{thm:rathomotopymaps}
For any $X$ and $Y$, $[X,Y]^G_\qq$ is a rational vector space.
If $Z$ is an $S^0_\mcal{\qq}$-local object of $G \mcal$ then
$Z$ has rational homotopy groups.
There is a natural isomorphism
$[X,Y]^G_\qq \cong [X \smashprod S^0_\mcal{\qq},Y \smashprod S^0_\mcal{\qq}]^G.$
\end{thm}
\begin{pf}
For each integer $n$ we have a self-map of $\cofrep S$
which represents multiplication by $n$ at the model category level,
applying $(-) \smashprod X$ we obtain a self-map of
$\cofrep S \smashprod X$. Since this map is an isomorphism
of rational homotopy groups it induces an isomorphism 
$n \co [X,Y]^G_\qq \to [X,Y]^G_\qq$, 
hence $[X,Y]^G_\qq$ is a rational vector space.
The homotopy groups of $Z$ can be given in terms of
$[\Sigma^p G/H_+, Z]^G$ for $p$ an integer and $H$
a subgroup of $G$. Since we have assumed that $Z$ is  
$S^0_\mcal{\qq}$-local, this homotopy group is isomorphic to 
$[\Sigma^p G/H_+, Z]^G_\qq$ which we now know is a rational
vector space.

The map $Y \smashprod S^0_\mcal{\qq} \to \fibrep_\qq (Y \smashprod S^0_\mcal{\qq})$
is a $\pi_*^\qq$-isomorphism between objects with
rational homotopy groups, hence it is a
$\pi_*$-isomorphism. For any $G$-spectrum $X$,
$X \smashprod S^0_\mcal{\qq}$ is rationally equivalent to $X$.
Combining these we obtain isomorphisms as below.
$$ \begin{array}{rcll}
[X,Y]^G_\qq
& \cong & [X \smashprod S^0{\qq},Y \smashprod S^0{\qq}]^G_\qq \\
& \cong & [X \smashprod S^0{\qq},\fibrep_\qq (Y \smashprod S^0{\qq})]^G \\
& \cong & [X \smashprod S^0{\qq},Y \smashprod S^0{\qq}]^G
\end{array} $$
\end{pf}

The following result gives a universal property for 
$S^0_\mcal \qq$. Note that if the map $f$ is a 
rational equivalence, then the lift in the proof below
is a rational equivalence between spectra with rational
homotopy groups and hence is a weak equivalence. 

\begin{lem}\label{lem:univpropsq}
Let $X$ be a spectrum with a map $f \co S \to X$ 
such that $\pi_n^H(X)$ is a rational vector space
for each subgroup $H$ and integer $n$. 
Then there is a map $S^0_\mcal \qq \to X$ in $\ho G \mcal$
such that the composite $S \to S^0_\mcal \qq \to X$
is equal to the map $f$ (in $\ho G \mcal$).
\end{lem}
\begin{pf}
By Theorem \ref{thm:rathomotopymaps} the map 
$\cofrep X \to \fibrep_\qq \cofrep X$ is a weak equivalence.
We then draw the diagram below and obtain a lifting 
$S^0_\mcal \qq \to \fibrep_\qq \cofrep X$ using the rational model structure
on $G \mcal$.
$$\xymatrix{
*+<0.5cm>{\cofrep S} \ar[r] \ar@{>->}[d]_{\sim_\qq} & 
\cofrep X \ar[r]^{\sim} & \fibrep_\qq \cofrep X \ar@{->>}[d] \\
S^0_\mcal \qq \ar[rr] & & {\ast} }$$
\end{pf}

\section{Splitting Rational Equivariant Spectra}\label{sec:eqstabhom}

We show how splittings of the category of rational equivariant spectra
correspond to idempotents of the rational Burnside ring. 
In particular, we know all such idempotents in the case of a finite group
and we have the idempotent $e_1$, constructed in 
Lemma \ref{lem:Geidemfamily}, which is in many cases a non-trivial idempotent.
For a compact Lie group $G$ the Burnside ring is 
defined to be $[S,S]^G$.
The following result is tom Dieck's isomorphism, 
see \cite[Chapter V, Lemma 2.10]{lms86} which references
\cite[Lemma 6]{tdieck77}. This result can be very useful 
when studying the Burnside ring of $G$.
Recall that 
$\fcal G$ is the set of subgroups of $G$
that have finite index in their normaliser. 
There is a topology on $\fcal G$ 
(induced by the Hausdorff metric on subsets of $G$)
such that the conjugation action of $G$ on $\fcal G$
is continuous, see \cite[Chapter V, Lemma 2.8]{lms86}. 

\begin{lem}\label{lem:tdisom}
Let $C( \fcal G / G, \qq)$ denote the ring
of continuous maps from the orbit space $\fcal G / G$ to $\qq$, 
where $\qq$ is
considered as a topological space with the discrete topology.
The map $[S,S]^G \to C( \fcal G / G, \qq)$
which takes $f$ to $(H) \mapsto \deg(f^H)$
induces an isomorphism of rings 
$[S,S]^G \otimes \qq \to C( \fcal G / G, \qq)$.
\end{lem}
In particular, for a finite group $G$, this
specifies an isomorphism 
$[S,S]^G \otimes \qq \to \prod_{(H) \leqslant G} \qq$.
Let $e_H \in [S,S]^G \otimes \qq$ be the idempotent corresponding to projection onto
factor $(H)$, then we have a finite orthogonal decomposition of $\id_S$
given by the collection $\{ e_H \}$ as $H$ runs over the conjugacy
classes of subgroups of $G$. We now give an isomorphism between the rational Burnside ring 
and self maps of $S$ in $\ho G \mcal_\qq$.

\begin{prop}
There is a ring isomorphism 
$[S,S]^G \otimes \qq \to [S,S]^G_\qq$
induced by $\id \co G \mcal \to G \mcal_\qq$. 
\end{prop}
\begin{pf}
The identity functor induces a ring map
$[S,S]^G \to [S,S]^G_\qq$ and since the right hand side 
is a rational vector space this induces
the desired map of rings. 
That this map is an isomorphism follows from the isomorphisms:
$[S, S]^G \otimes \qq \cong [S, S^0_\mcal \qq]^G$,
$[S, S^0_\mcal \qq]^G \cong [S, \fibrep_\qq S]^G$ and
$[S, \fibrep_\qq S]^G \cong [S, S]^G_\qq$. 
The universal property of $S^0_\mcal \qq$ provides 
the second isomorphism and ensures that the composite of the above maps
is equal to the specified map of rings.  
\end{pf}

\begin{cor}\label{cor:Gspecsplit}
If $e$ is an idempotent of the rational Burnside ring of $G$, 
then the adjunction below is a strong symmetric monoidal Quillen
equivalence. 
$$\Delta :  G \mcal_\qq \overrightarrow{\longleftarrow} 
L_{e S} G \mcal_\qq  \times L_{(1-e) S} G \mcal_\qq  : \prod$$
\end{cor}

\begin{cor}\label{cor:finG}
The category of rational $G$-spectra 
(for finite $G$) splits into the product of the localisations
$L_{e_H S} \GIS_\qq$ as $(H)$ runs over the conjugacy
classes of subgroups of $G$. 
\end{cor}
At the homotopy level 
this result can be found in \cite[Appendix A]{gremay95}.
Note that the two localisations of $G$-spectra that we have used:
$L_{S^0_\mcal \qq} G \mcal$ and $L_{eS} L_{S^0_\mcal \qq} G \mcal$ share many of the same
properties. This is because they are designed to invert elements
of $[S,S]^G$ and $[S,S]_\qq^G$ respectively. 
The first is designed to invert 
the primes and the second 
inverts the idempotent $e$. 

\begin{lem}\label{lem:rightproperfamily}
For $e$ an idempotent of $[S,S]^G \otimes \qq$ 
the category $L_{eS} G \mcal$ is right proper. 
\end{lem}
\begin{pf}
Let $e \in [S,S]^G \otimes \qq$ be an idempotent, then
for any exact sequence of $[S,S]^G \otimes \qq$-modules
$ \dots \to M_i \to M_{i-1} \to \dots $, the sequence
$ \dots \to e M_i \to e M_{i-1} \to \dots $ is exact.
Right properness then follows from the proof of
Lemma \ref{lem:rightproperrational} by
applying $e$ to the long exact sequence 
of rational homotopy groups
of a fibration. 
\end{pf}

We now give a general example of an idempotent of 
the Burnside ring. This idempotent is non-trivial
in many cases, such as when $G=O(2)$, the group of 
two-by-two orthogonal matrices. This idempotent was
used to study rational $O(2)$-spectra in 
\cite{gre98a} and \cite[Part III]{barnes}.

\begin{lem}\label{lem:Geidemfamily}
Let $G$ be a compact Lie group
and let $S$ denote the set of subgroups
of the identity component of $G$
which have finite index in their normaliser. 
Then there is an idempotent $e_1 \in [S,S]^G \otimes \qq \cong C( \fcal G / G, \qq)$
given by the map which sends $(H)$ to 1
if $H \in S$ and zero otherwise.
\end{lem}
\begin{pf}
Let $G_1$ denote the identity component of $G$ and
recall that since $G$ is compact
$F = G /G_1$ is finite.
Take $H \in S$, by \cite[Chapter II, Corollary 5.6]{bred} we know that
if $K \in \fcal G$ is in some sufficiently small neighbourhood of
$H$ in the space $\fcal G$, then $K$ is subconjugate
to $H$ and so $K$ is a subgroup of $G_1$. It follows that
$S$ is open in $\fcal G/G$.
Now take $(K)$ to be in $(\fcal G/G) \setminus S$,
so there is a $g \in G \setminus G_1$ such that $K \cap gG_1$ is non-empty.
Then any $L \in \fcal G$ that is sufficiently close to $K$
also has a non-trivial intersection with $gG_1$
so $L$ is not a subgroup of $G_1$, 
it follows that $S$ is also closed.
Hence $e_1$, 
the characteristic function of $S$,
is a continuous map $\fcal G /G \to \qq$.
Thus $e_1$ is an idempotent, since $e_1(H)=1$
if $H \in S$ and zero otherwise. 
\end{pf}

Let $\fscr$ be the set of subgroups
of $G_1$, then it can be shown that $e_1 S$ is weakly equivalent
to ${E \fscr}_+$ (the universal space for a family). 
One can then use the results of 
\cite[Chapter IV, Section 6]{mm02} to obtain
better understanding of $L_{e_1 S} G \mcal_\qq$
and $L_{(1-e_1)S} G \mcal_\qq$.

\section{Modules and Bimodules}\label{sec:modbimod}
We give two general examples of where our splitting result can be applied. 
Choose a monoidal model category of spectra, such as
symmetric, orthogonal or EKMM spectra (this could even 
be $G$-equivariant for the last two versions) and call it 
$\sscr$. For $R$ a ring spectrum we consider 
splittings of the model category 
of $R$-$R$-bimodules, this is a monoidal
model category which is not (in general) symmetric. 
We let $[-,-]^{(R,R)}$ denote maps in the homotopy category of $R$-$R$-bimodules.
Our second example considers the case of $R$-modules, when $R$
is not commutative. Although $R \leftmod$ is not
a monoidal model category we can still 
obtain splittings of the model category
by considering idempotents of $[R,R]^{(R,R)}$. 
We return to rational equivariant spectra
at the end of this section and create a commutative ring spectrum
$S_\qq$ such that $S_\qq \leftmod$ is Quillen equivalent to 
$G \mcal_\qq$ (Theorem \ref{thm:localisedtomodules}). 
We then show that splittings of 
$S_\qq \leftmod$ correspond to splittings of 
$G \mcal_\qq$.

We first introduce some results from 
\cite{EKMM97}, these can be adapted to any of the categories
of spectra we have mentioned above. 
For $R$ an algebra, there is a notion of a cell $R$-module, 
see \cite[Chapter III, Definition 2.1]{EKMM97}, a cell $R$ module
is a special kind of cofibrant module. 
We can always replace an $R$-module $M$
by a weakly equivalent cell $R$-module $\Gamma M$
via \cite[Chapter III, Theorem 2.10]{EKMM97}.

If $E$ is a right $R$-module 
then we have a spectrum $E \smashprod_R X$
for any left $R$-module $X$. It is defined as the 
coequaliser of the diagram 
$E \smashprod R \smashprod X \overrightarrow{\longrightarrow} E \smashprod X$
where the maps are given by the action of $R$ on $E$ 
and the action of $R$ on $X$. 
Thus we have the notion of an 
$E^R$-equivalence of $R$-modules: a map $f$ in 
$R \leftmod$ such that $E \smashprod_R f$ is a weak equivalence
of underlying spectra. 
Let $E$ be a cell right $R$-module, then by 
\cite[Chapter VIII, Theorem 1.1]{EKMM97}, 
there is a model structure 
$L_E R \leftmod$ on the category of $R$-modules
with weak equivalences the $E^R$-equivalences
and cofibrations given by the cofibrations for $R \leftmod$.
We also note that if $X$ is a cofibrant $R$-module,
the functor $- \smashprod_R X$ preserves weak equivalences 
(\cite[Chapter III, Theorem 3.8]{EKMM97}).

\begin{prop}\label{prop:splitbimod}
For $R$ a ring spectrum in $\sscr$, 
whose underlying spectrum is cofibrant,
an idempotent of $\text{THH}^0 (R):=[R,R]^{(R,R)}$
splits the category of $R$-$R$-bimodules. 
\end{prop}
\begin{pf}
We can identify the category of $R$-$R$ bimodules
with the category of $R \smashprod R^{op}$-modules.
The ring spectrum $R^{op}$ has the same underlying spectrum as $R$
but the multiplication is given by 
$R \smashprod R \overset{\tau}{\to} R \smashprod R \overset{\mu}\to R$
where $\tau$ is the symmetry isomorphism of $\smashprod$ in $\sscr$
and $\mu$ is the multiplication of $R$. 
We have assumed that $R$ is cofibrant to ensure that 
$R \smashprod R^{op}$ is weakly equivalent to 
$R \smashprod^L R^{op}$, thus  
$[X,Y]^{(R,R)} \cong [X,Y]^{R \smashprod^L R^{op}}$.

For a cell $R$-$R$-bimodule $E$ 
we have a $E$-local model structure on the category of $R$-$R$-bimodules.
If $M$ is a cofibrant $R$-$R$-bimodule, then an $M$-equivalence is the same
as a $\Gamma M$-equivalence 
and so we can localise at any cofibrant bimodule
by localising at its cellular replacement.
We can now apply Theorem \ref{thm:generalsplitting} 
to complete the proof. 
\end{pf}

We now turn to left modules over a ring spectrum, 
we can obtain a splitting result
when $R$ is not commutative. In this case $R \leftmod$ 
does not have a monoidal product and so
$[R,R]^R$ does not act on $[X,Y]^R$. 
Instead we will use the action of $[R,R]^{(R,R)}$ on 
$[X,Y]^R$ to split the category. Throughout
we assume that $R$ is cofibrant as a spectrum. 

We return to algebra briefly to offer some context for this result. 
If $R$ was an arbitrary ring, then for a \emph{central} idempotent
$e \in R$, (so $er=re$ for any $r \in R$), one can form
new rings $eR$ and $(1-e)R$ such that 
$R \leftmod$ is equivalent to $eR \leftmod \times (1-e)R \leftmod$. 
Furthermore, for any $R$-module $M$, there is a natural isomorphism 
$M \cong eM \oplus (1-e)M$. A central idempotent
is precisely the same data as an $R$-$R$-bimodule map from 
$R$ to itself. Hence, the proposition below is the ring spectrum
version of this algebraic result. 

\begin{prop}\label{prop:splitrmod}
Let $R \in \sscr$ be a ring spectrum whose underlying spectrum is cofibrant
and let $e$ be an idempotent of 
$[R,R]^{(R,R)}$. Then there is a Quillen equivalence
$$\Delta : R \leftmod \overrightarrow{\longleftarrow}
L_{\Gamma e R} R \leftmod \times L_{\Gamma (1-e) R} R \leftmod : \prod. $$
\end{prop}
\begin{pf}
We construct $e R$ in the category of $R$-$R$-bimodules
and then consider it as a right $R$-module. 
Since $R$ is cofibrant, it follows that $eR$ is cofibrant
as a right $R$-module (see below for details). 
We localise the category of $R$-modules at the cell right $R$-module
$\Gamma eR$ and note that the weak equivalences of 
$L_{\Gamma e R} R \leftmod$ are the $(eR)^R$-equivalences.
We can then follow the proof of Theorem \ref{thm:generalsplitting}.
\end{pf}

There is a forgetful functor $U$
from $R$-$R$-bimodules to $R \leftmod$, this is a right
Quillen functor with left adjoint
$M \mapsto M \smashprod R$.
Take $f: A \to B$ a generating (acyclic) cofibration
of $\sscr$. Then $g= \id_R \smashprod f \smashprod \id_R$
is a generating (acyclic) cofibration for the category of 
$R$-$R$-bimodules. Since $f \smashprod \id_R$
is a cofibration of spectra, it follows that 
$g$ is a cofibration of left $R$-modules, hence
$U$ is a left Quillen functor. 
A slight alteration of this argument shows that 
a cofibrant $R$-$R$-bimodule is 
cofibrant as a right $R$-module. 

The functor $U$ induces a ring map 
$[R,R]^{(R,R)} \to [R,R]^R \cong \pi_0(R)$.
If $R$ is commutative, every $R$-module
can be considered as an $R$-$R$-bimodule, 
this defines a right Quillen functor $I$.
Let $M$ be an $R$-$R$-bimodule with actions 
$\nu$ and $\nu'$. Then define $SM$ as the coequaliser:
$\xymatrix@C+0.3cm{
R \smashprod M 
\ar@<+0.2ex>[r]^(0.6){\nu}
\ar@<-0.2ex>[r]_(0.6){\nu' \circ \tau} &
M \ar[r] &
SM. }$
It follows that $S$ is the left adjoint of $I$
and that $UI$ is the identity functor of 
$R \leftmod$. 

These functors give a retraction:
$[R,R]^R \overset{I}{\to}
[R,R]^{(R,R)} \overset{U}{\to}
[R,R]^R$. Thus in the commutative case
it is no restriction to consider an
idempotent $e \in [R,R]^{(R,R)}$.
The Quillen equivalence above would then follow from our
main result and would be a 
strong symmetric monoidal Quillen equivalence.

For $E$ a cofibrant spectrum and $R$ a commutative ring spectrum, 
the $L_{E \smashprod R}$-model 
structure on the category of $R$-modules
has weak equivalences those maps $f$
which are $E$-equivalences of underlying spectra.
Thus $L_{E \smashprod R} R \leftmod$
is precisely the model category of $R$-modules
in $L_E \sscr$.

One important source of idempotents 
in $\pi_0 (R)$ (or $[R,R]^{(R,R)}$) is the image of idempotents
in $\pi_0(S)$ via the unit map $S \to R$. 
We return to our primary example of rational equivariant EKMM-spectra
to give an example of this.
To obtain our commutative ring spectrum
we use \cite[Chapter VIII, Theorem 2.2]{EKMM97},
we give the statement that we will need
below. Here we assume that $E$ is
a cell spectrum (hence cofibrant).

\begin{thm}\label{thm:algebralocalise}
For a cell commutative $R$-algebra $A$, the localisation
$\lambda : A \to A_E$ can be constructed as the inclusion of
a subcomplex in a cell commutative $R$-algebra $A_E$.
In particular $A \to A_E$ is an $E$-equivalence
and a cofibration of commutative ring spectra
for any cell commutative $R$-algebra $A$.
\end{thm}

\begin{defn}\label{def:S_qq}
Let $S_\qq$\index{Sa @$S_\qq$} be the commutative ring spectrum
constructed as the $S^0_\mcal{\qq}$-localisation of $S$.
\end{defn}
It follows immediately that the unit
$\eta : S \to S_\qq$ is an $S^0_\mcal \qq$-equivalence.
Thus, by our universal property for $S^0_\mcal \qq$
(Lemma \ref{lem:univpropsq})
and the fact that $S_\qq$ has rational homotopy groups,
we have the first statement of the following result. 
The rest of the lemma follows by a standard argument, see
\cite[13.1]{adams}. 

\begin{lem}\label{lem:compareSQtoRQ}
There is a weak equivalence
$S^0_\mcal{\qq} \to S_\qq$.
Hence all $S_\qq$-modules are $S^0_\mcal{\qq}$-local and 
so all $S_\qq$-modules have rational homotopy groups.
\end{lem}

\begin{thm}\label{thm:localisedtomodules}
There is a strong symmetric monoidal Quillen equivalence:
$$S_\qq \smashprod (-) :  G \mcal_\qq 
\overrightarrow{\longleftarrow} S_\qq \leftmod : U.$$
\end{thm}
\begin{pf}
The above functors form a strong monoidal Quillen pair
(with the usual structure on $G \mcal $).
Since cofibrations are unaffected
by localisation, $S_\qq \smashprod (-) : G \mcal_\qq \to S_\qq \leftmod$ preserves
cofibrations. Consider an acyclic rational cofibration $X \to Y$,
we know that $S_\qq \smashprod (-)$ applied
to this gives a cofibration, we must check that
it is also a $\pi_*$-isomorphism. 

We see that
$X \smashprod S^0_\mcal{\qq} \to Y \smashprod S^0_\mcal{\qq}$ is a
cofibration and a $\pi_*$-isomorphism,
so in turn
$X \smashprod S^0_\mcal{\qq} \smashprod S_\qq \to
Y \smashprod S^0_\mcal{\qq} \smashprod S_\qq$ is 
a $\pi_*$-isomorphism (by the monoid axiom).
This proves that
$X \smashprod S_\qq \to Y \smashprod S_\qq$
is a $\pi^\qq_*$-isomorphism between
$S_\qq$-modules, which we know have rational homotopy groups and thus
this map is in fact a $\pi_*$-isomorphism. Hence we have a Quillen pair, 
now we prove that it is a Quillen equivalence.
The right adjoint preserves and detects all weak equivalences.
The map $X \to S_\qq \smashprod X$ is a rational equivalence 
for all cofibrant $S$-modules $X$.
This follows since smashing with a cofibrant 
object will preserve the $\pi_*^\qq$-isomorphism $S \to S_\qq$. 
\end{pf}

It follows that we have an isomorphism of rings 
$[S,S]^G_\qq \to [S_\qq, S_\qq]^{S_\qq \leftmod}$. 
Hence for an idempotent $e$ of the rational Burnside ring 
we can split $S_\qq \leftmod$ using
the objects $e S \smashprod S_\qq$ and 
$(1-e)S \smashprod S_\qq$. 
We can then apply Theorem \ref{thm:locfuncs} to see that
the strong symmetric monoidal adjunction below is a Quillen equivalence, 
hence we have a comparison between our splitting of $S_\qq \leftmod$
and Corollary \ref{cor:Gspecsplit}.
$$S_\qq \smashprod (-): L_{eS} G \mcal_\qq
\overrightarrow{\longleftarrow}
L_{(\Gamma eS) \smashprod S_\qq} S_\qq \leftmod :U$$
We briefly wish to mention that following the construction of $S_\qq$
one can make $R_{e}$ for any commutative ring $R$ and idempotent
$e \in \pi_0(R)$ by localising $R$ at $\Gamma eR$. 
It follows that $R_e$ is weakly equivalent to $\Gamma eR$
and hence any $R_e$-module is $\Gamma eR$-local. 
Then, as with the $S_\qq$-case, one can prove that 
extension and restriction of scalars along 
$R \to R_e$ induces a Quillen equivalence between
$L_{\Gamma eR} R \leftmod$ and $R_e \leftmod$. 
This is a manifestation of \cite[Theorem 2]{wolb}.
Hence we have a different statement of the splitting result:
there is a Quillen equivalence $R \leftmod 
\overrightarrow{\longleftarrow}
R_e \leftmod \times R_{1-e} \leftmod$, 
induced by extension and restriction of
scalars.

\bibliography{daveprebib}
\bibliographystyle{alpha}

\end{document}